\documentclass[11pt,a4paper]{article}
\usepackage{mathrsfs}
\usepackage{amsmath}
\usepackage{amsfonts}
\usepackage{amssymb}

\usepackage{graphicx,subfigure,epsfig}
\usepackage{enumerate}
\usepackage{arydshln}
\usepackage[pdftex,linkcolor=blue,citecolor=blue,backref=pagebackref]{hyperref}
\usepackage{color}
\usepackage{pdfpages}
\usepackage{subfigure}
\usepackage{tabularx}

\allowdisplaybreaks[4]

\begin{document}

\renewcommand{\thefootnote}{\fnsymbol{footnote}}

\newtheorem{theorem}{Theorem}[section]
\newtheorem{corollary}[theorem]{Corollary}
\newtheorem{definition}[theorem]{Definition}
\newtheorem{conjecture}[theorem]{Conjecture}
\newtheorem{question}[theorem]{Question}
\newtheorem{lemma}[theorem]{Lemma}
\newtheorem{proposition}[theorem]{Proposition}
\newtheorem{example}[theorem]{Example}
\newtheorem{fact}[theorem]{Fact}
\newenvironment{proof}{\noindent {\bf
Proof.}}{\rule{3mm}{3mm}\par\medskip}
\newcommand{\remark}{\medskip\par\noindent {\bf Remark.~~}}
\newcommand{\pp}{{\it p.}}
\newcommand{\de}{\em}
\newcommand{\qbinom}[2]{\genfrac{[}{]}{0pt}{}{#1}{#2}}

\newcommand{\JEC}{{\it Europ. J. Combinatorics},  }
\newcommand{\JCTB}{{\it J. Combin. Theory Ser. B.}, }
\newcommand{\JCT}{{\it J. Combin. Theory}, }
\newcommand{\JGT}{{\it J. Graph Theory}, }
\newcommand{\ComHung}{{\it Combinatorica}, }
\newcommand{\DM}{{\it Discrete Math.}, }
\newcommand{\ARS}{{\it Ars Combin.}, }
\newcommand{\SIAMDM}{{\it SIAM J. Discrete Math.}, }
\newcommand{\SIAMADM}{{\it SIAM J. Algebraic Discrete Methods}, }
\newcommand{\SIAMC}{{\it SIAM J. Comput.}, }
\newcommand{\ConAMS}{{\it Contemp. Math. AMS}, }
\newcommand{\TransAMS}{{\it Trans. Amer. Math. Soc.}, }
\newcommand{\AnDM}{{\it Ann. Discrete Math.}, }
\newcommand{\NBS}{{\it J. Res. Nat. Bur. Standards} {\rm B}, }
\newcommand{\ConNum}{{\it Congr. Numer.}, }
\newcommand{\CJM}{{\it Canad. J. Math.}, }
\newcommand{\JLMS}{{\it J. London Math. Soc.}, }
\newcommand{\PLMS}{{\it Proc. London Math. Soc.}, }
\newcommand{\PAMS}{{\it Proc. Amer. Math. Soc.}, }
\newcommand{\JCMCC}{{\it J. Combin. Math. Combin. Comput.}, }
\newcommand{\GC}{{\it Graphs Combin.}, }
\newcommand{\LAA}{{\it Linear Algeb. Appli.}, }

\title{ \bf A Relationship for LYM Inequalities between
\\ Boolean Lattices and Linear Lattices with Applications}

\author{Jiuqiang Liu$^{a,b,} \thanks{The corresponding author}$, Guihai Yu$^{a, *}$\\
{\small  $^a$ College of Big Data Statistics, Guizhou University of Finance and Economics}\\
 {\small Guiyang, Guizhou, 550025, China}\\
 {\small $^{b}$ Department of Mathematics, Eastern Michigan University}\\
{\small Ypsilanti, MI 48197, USA}\\
{\small E-mail: { \tt jiuqiang68@126.com, yuguihai@mail.gufe.edu.cn}}}

\maketitle
\vspace{-0.5cm}

\begin{abstract}
Sperner theory is one of the most important branches in extremal set theory. It has many applications in the field of operation research, computer science, hypergraph theory and so on. The LYM property has become an important tool for studying Sperner property.
In this paper, we provide a general relationship for LYM inequalities between Boolean lattices and linear lattices. As applications, we use this relationship to derive generalizations of some well-known theorems on maximum sizes of families containing no copy of certain poset or certain configuration from Boolean lattices to linear lattices, including generalizations of the well-known Kleitman theorem on families containing no $s$ pairwise disjoint members (a non-uniform variant of the famous Erd\H{o}s matching conjecture)
and Johnston-Lu-Milans theorem and Polymath theorem on families containing no $d$-dimensional Boolean algebras.
\end{abstract}

 \vskip 0.1cm

\section{Introduction}
\hspace*{0.5cm}
Throughout this paper, we denote $[n] = \{1, 2, \dots, n\}$. A {\em Boolean lattice} $\mathcal{B}_{n}$ is the set $2^{[n]}$ of all subsets of $[n]$ with the ordering being the containment relation and a {\em linear lattice}
$\mathcal{L}_{n}(q)$ is the set of all subspaces of the $n$-dimensional vector space $\mathbb{F}_{q}^{n}$ over the finite field $\mathbb{F}_{q}$ with the ordering being the relation of inclusion of subspaces.
Denote by ${{[n]} \choose {k}}$ for the set of all $k$-subsets of $[n]$ and $\qbinom{[n]}{i}$ denotes the set of all $i$-dimensional subspaces of $\mathbb{F}_{q}^{n}$.

Given posets $P$ and $P'$, we say that $P$ is {\em weakly contained} in $P'$ if there exists an injective map
$\psi : P \longrightarrow P'$ such that for every $u, v \in P$, $\psi(u) \leq _{P'} \psi(v)$ if $u \leq_{P} v$; and $P'$ {\em strongly contains} $P$ if
for all $u, v \in P$, $\psi(u) \leq _{P'} \psi(v)$ if and only if $u \leq_{P} v$.

Given a poset $P$, $ex(n, P)$ is defined by
\[ex(n, P) = \max \{|\mathcal{F}| \mid \mathcal{F} \subseteq 2^{[n]} \mbox{ and } \mathcal{F} \mbox{ does not contain weakly } P\}\]
and $en^{*}(n, P)$ is defined by
\[ex^{*}(n, P) = \max \{|\mathcal{F}| \mid \mathcal{F} \subseteq 2^{[n]} \mbox{ and } \mathcal{F} \mbox{ does not contain strongly } P.\}\]
There has been a considerable amount of research devoted to determining the asymptotic
behaviour of $ex(n, P)$ and $ex^{*}(n, P)$. The first result of this kind is the following well-known Sperner's Theorem.
Recall that a family of subsets of $[n]$ is called a {\em Sperner family} (or antichain) if there
are no two different members of the family such that one of them contains the
other.

\begin{theorem}\label{thm1.1} (Sperner Theorem). Suppose that $\mathcal{A}$ is a Sperner family
of subsets of $[n]$. Then
\[|\mathcal{A}|\leq {{n} \choose {\lfloor \frac{n}{2}\rfloor}}.\]
\end{theorem}
Clearly,, Sperner's result is equivalent to say that $ex(n, P_{2}) =  {{n} \choose {\lfloor \frac{n}{2}\rfloor}}$.

We say a family $\mathcal{F}$ of subsets of $[n]$ is $k$-$Sperner$
if all chains in $\mathcal{F}$ have length at most $k$. Define $\sum (n, k)$ to be the sum of the
$k$ largest binomial coefficients of order $n$, i.e., $\sum (n, k) = \sum_{i = 1}^{k}{{n} \choose {\lfloor \frac{n - k}{2}\rfloor + i}}$.
Let $\sum^{*} (n, k)$ be the collection of families consisting of the corresponding full levels,
i.e., if $n + k$ is odd, then $\sum^{*} (n, k)$ contains one family $\cup_{i = 1}^{k}{{[n]} \choose {\lfloor \frac{n - k}{2}\rfloor + i}}$;
if $n + k$ is even, then $\sum^{*} (n, k)$ contains two families of the same size $\cup_{i = 0}^{k - 1}{{[n]} \choose {\frac{n - k}{2}+ i}}$ and
$\cup_{i = 1}^{k}{{[n]} \choose {\frac{n - k}{2}+ i}}$. The following theorem by Erd\H{o}s \cite{e} generalizes the classical Sperner theorem.

\begin{theorem}\label{thm1.2} (Erd\H{o}s, \cite{e}). Suppose that $\mathcal{A}$ is a $k$-Sperner family
of subsets of $[n]$. Then
\[|\mathcal{A}|\leq \sum (n, k).\]
Moreover, if $|\mathcal{A}| = \sum (n, k)$, then $\mathcal{A} \in \sum^{*} (n, k)$.
\end{theorem}

Denote
\begin{equation}
\qbinom{n}{k} = \prod_{0 \leq i \leq k - 1}\frac{q^{n - i} - 1}{q^{k - i} - 1}.
\label{1.1}\end{equation}
It is well known that the number of all $k$-dimensional
subspaces of $\mathbb{F}_{q}^{n}$ is equal to $\qbinom{n}{k}$.

Note that there has been a considerable amount of research on maximum sizes of families of subsets of $[n]$ forbidden certain structures and many of the results appeared in the literature have their $q$-analogues (generalizations to linear lattices) which are usually derived using different and more complicated methods, see for example \cite{abs}, \cite{cp}, \cite{ekr}, \cite{fg}, \cite{fw}, \cite{gk}, \cite{h}, \cite{qr}, and \cite{rw}.

We say that a family $\mathcal{V}$ of subspaces is {\em Sperner} (or antichain) if no subspace is contained in another subspace in $\mathcal{V}$, and $\mathcal{V}$ is
$k$-{\em Sperner} if all chains in $\mathcal{V}$ have length at most $k$.
The following well-known Sperner theorem for families of subspaces of $\mathbb{F}_{q}^{n}$ can be found in \cite{e2}, which is a vector space analogue (or $q$-analogue) of the classical Sperner theorem.

\begin{theorem}\label{thm1.3} (q-Analogue Sperner Theorem).  Assume that $\mathcal{V}$ is a Sperner family
of subspaces of $\mathbb{F}_{q}^{n}$. Then
\[|\mathcal{V}| \leq \qbinom{n}{\lfloor \frac{n}{2}\rfloor}.\]
\end{theorem}

We denote $\sum [n, k]$ to be the sum of the
$k$ largest $q$-binomial coefficients of order $n$, i.e., $\sum [n, k] = \sum_{i = 1}^{k} \left[n \atop \lfloor \frac{n - k}{2}\rfloor + i \right]$.
Let $\sum^{*} [n, k]$ be the collection of families consisting of the corresponding full levels,
i.e., if $n + k$ is odd, then $\sum^{*} [n, k]$ contains one family $\cup_{i = 1}^{k}\left[[n] \atop \lfloor \frac{n - k}{2}\rfloor + i \right]$;
if $n + k$ is even, then $\sum^{*} [n, k]$ contains two families of the same size $\cup_{i = 0}^{k - 1}\left[[n] \atop \frac{n - k}{2} + i \right]$
and $\cup_{i = 1}^{k}\left[[n] \atop \frac{n - k}{2} + i \right]$.

The following $q$-analogue of Erd\H{o}s' theorem (Theorem \ref{thm1.2}) is implied by Theorem 2 in \cite{s2}, which generalizes Theorem \ref{thm1.3} to $k$-Sperner families.

\begin{theorem}\label{thm1.4} (Samotij, \cite{s2}).  Assume that $\mathcal{V}$ is a $k$-Sperner family
of subspaces of $\mathbb{F}_{q}^{n}$. Then
\[|\mathcal{V}|\leq \sum [n, k].\]
Moreover, if $|\mathcal{V}| = \sum [n, k]$, then $\mathcal{V} \in \sum^{*} [n, k]$.
\end{theorem}

Given $\mathcal{F} \subseteq 2^{[n]}$, the {\em Lubell function} of $\mathcal{F}$, also called {\em weight} and denoted by $l(\mathcal{F})$, is the quantity
\begin{equation}
l(\mathcal{F}) = \sum_{F \in \mathcal{F}}\frac{1}{{{n} \choose {|F|}}} = \sum_{i = 0}^{n}\frac{|\mathcal{F}_{i}|}{{{n} \choose {i}}},
\end{equation}
where, $\mathcal{F}_{i} = \{ F \in \mathcal{F} \mid |F| = i\}$.

One of the advantages of using the Lubell function is its convenient probabilistic interpretation: Suppose that
$\mathcal{C} = \{\emptyset, \{i_{1}\}, \{i_{1}, i_{2}\}, \dots, [n]\}$
is full-chain in $2^{[n]}$ chosen uniformly at random.
Let $X$ be the random variable $X = |\mathcal{C} \cap \mathcal{F}|$. Then one has the expected value $E(X) = l(\mathcal{F})$.

The following well-known LYM inequality was proved independently by Bollob\'{a}s \cite{b}, Lubell \cite{l2}, Meshalkin \cite{m1} and
Yamamoto \cite{y}.

\begin{theorem}\label{thm1.5} (LYM Inequality, Bollob\'{a}s \cite{b}, Lubell \cite{l2}, Meshalkin \cite{m1} and
Yamamoto \cite{y}).
If $\mathcal{F} \subseteq 2^{[n]}$ is an antichain, then
\[l(\mathcal{F}) = \sum_{F \in \mathcal{F}}\frac{1}{{{n} \choose {|F|}}} \leq 1.\]
\end{theorem}

This LYM inequality naturally implies Sperner's Theorem (Theorem \ref{thm1.1}) and is a cornerstone
of the field of extremal set theory, and the weight or Lubell function $l(\mathcal{F})$ is a powerful tool for studying both $ex(n, P)$ and $ex^{*}(n, P)$ and has been used in various papers, see for example \cite{gll}, \cite{gmt}, \cite{jlm}, and \cite{kmy}. The LYM inequality has been strengthened by Bey \cite{b1} and  Erd\H{o}s et al. \cite{efkss} and has a continuous analogue by Klain and Rota \cite{kr}.

In general, given a ranked poset $\mathcal{P}$, in which the set of
elements of rank $r$ is denoted by $\mathcal{P}_{r}$, $\mathcal{P}$ has the {\em Sperner property} if the largest antichain has size equal to the size of the largest rank level, and $\mathcal{P}$ satisfies the {\em LYM property} (or LYM inequality) if for every antichain $\mathcal{A}$ of $\mathcal{P}$,
\begin{equation}
\sum_{i}\frac{|\mathcal{A}\cap \mathcal{P}_{i}|}{|\mathcal{P}_{i}|} \leq 1.
\end{equation}
Clearly, for $\mathcal{P} = \mathcal{B}_{n}$, the inequality (3) is the same as the inequality in Theorem \ref{thm1.5} with $\mathcal{P}_{i} = {{[n]} \choose {i}}$.

For $k$-Sperner families, one can show the next fact.

\begin{theorem}\label{thm1.6}
If $\mathcal{F} \subseteq 2^{[n]}$ contains no $P_{k+1}$, then
\[l(\mathcal{F}) \leq k.\]
\end{theorem}

This inequality implies Erd\H{o}s $k$-Sperner theorem (Theorem \ref{thm1.2}).
In 2018, M$\grave{e}$roueh \cite{m} proved the following general result which was conjectured to be true by Lu and Milans \cite{lm}.

\begin{theorem}\label{thm1.7} (M$\grave{e}$roueh, \cite{m}).
For every poset $P$, there exists $c(P)$ such that if $\mathcal{F} \subseteq 2^{[n]}$ for some $n \in \mathbb{N}$
and $\mathcal{F}$ does not contain $P$ as an induced subposet, then $l(\mathcal{F}) \leq c(P)$.
\end{theorem}

Clearly, a family $\mathcal{F} \subseteq 2^{[n]}$ does not contain $P$ as a subposet implies it does not contain $P$ as an induced subposet.
Theorem \ref{thm1.7} has the next immediate consequence.

\begin{theorem}\label{thm1.8}
For every poset $P$, there exists $c(P)$ such that if $\mathcal{F} \subseteq 2^{[n]}$ for some $n \in \mathbb{N}$
and $\mathcal{F}$ does not contain $P$ as a subposet, then $l(\mathcal{F}) \leq c(P)$.
\end{theorem}

Our first main result is the following general relationship for LYM inequalities between Boolean lattices and linear lattices (many results about maximum sizes of families of subspaces of vector space $\mathbb{F}_{q}^{n}$ can be derived from the corresponding theorems for maximum sizes of families of subsets of $[n]$ through this relationship).
For a family $\mathcal{V}$ of subspaces of $\mathbb{F}_{q}^{n}$, we denote
\begin{equation}
l_{q}(\mathcal{V}) = \sum_{V \in \mathcal{V}}\frac{1}{\qbinom{n}{|V|}} = \sum_{i = 0}^{n}\frac{|\mathcal{V}_{i}|}{\qbinom{n}{i}},
\end{equation}
where $\mathcal{V}_{i} = \{ V \in \mathcal{V} \mid dim(V) = i\}$.

\begin{theorem}\label{thm1.9}
Let $P$ be a given poset or configuration. Suppose that there exists $c(P)$ such that if $\mathcal{F} \subseteq 2^{[n]}$
and $\mathcal{F}$ does not contain $P$, then $l(\mathcal{F}) \leq c(P)$.
Then every family $\mathcal{V}$ of subspaces of $\mathbb{F}_{q}^{n}$
which contains no copy of $P$ satisfies $l_{q}(\mathcal{V}) = \sum_{V \in \mathcal{V}}\frac{1}{\qbinom{n}{|V|}} \leq c(P)$.
\end{theorem}

Clearly, Theorems \ref{thm1.8} and \ref{thm1.9} imply the next consequence which is a generalization of Theorem \ref{thm1.8} to linear lattices.

\begin{theorem}\label{thm1.10}
For every poset $P$, there exists $c(P)$ such that
every family $\mathcal{V}$ of subspaces of $\mathbb{F}_{q}^{n}$
which contains no copy of $P$ satisfies $l_{q}(\mathcal{V}) = \sum_{V \in \mathcal{V}}\frac{1}{\qbinom{n}{|V|}} \leq c(P)$.
\end{theorem}

Define the weighted Lubell functions $l^{*}(\mathcal{F})$ and $l^{*}_{q}(\mathcal{V})$ by
\[l^{*}(\mathcal{F}, \beta) = \sum_{i = 0}^{n}\frac{\beta_{i}|\mathcal{F}_{i}|}{{{n} \choose {i}}},\]
\[l^{*}_{q}(\mathcal{V}, \beta) = \sum_{i = 0}^{n}\frac{\beta_{i}|\mathcal{V}_{i}|}{\qbinom{n}{i}},\]
where $\mathcal{F}_{i} = \mathcal{F} \cap {{[n]} \choose {i}}$ and $\mathcal{V}_{i} = \mathcal{V} \cap \qbinom{[n]}{i}$.
We have the next weighted version of Theorem \ref{thm1.9}.

\begin{theorem}\label{thm1.11}
Let $P$ be a given poset or configuration. Suppose that there exists $c(P)$ such that if $\mathcal{F} \subseteq 2^{[n]}$
and $\mathcal{F}$ does not contain $P$, then $l^{*}(\mathcal{F}, \beta) \leq c(P)$.
Then every family $\mathcal{V}$ of subspaces of $\mathbb{F}_{q}^{n}$
which contains no copy of $P$ satisfies $l^{*}_{q}(\mathcal{V}, \beta) \leq c(P)$.
\end{theorem}

As applications of Theorem \ref{thm1.9}, in sections 3 - 5, we derive generalizations of some well-known theorems on maximum sizes of families containing no copy of certain poset or certain configuration from Boolean lattices to linear lattices, including generalizations of Gr\'osz-Methuku-Tompkins theorem on diamond-free families, generalizations of Johnston-Lu-Milans theorem and Polymath theorem on families containing no $d$-dimensional Boolean algebras which is used by Polymath in a new proof of  Furstenberg and Katznelson's density Hales-Jewett theorem (a generalization of the well-known Szemer\'edi theorem), and
generalizations of the following well-known Kleitman theorem on families containing no $s$ pairwise disjoint members (a non-uniform variant of the well-known Erd\H{o}s Matching Conjecture \cite{e1}) and Frankl-Kupavskii theorem on cross-dependent families \cite{fk} as shown in Theorems \ref{thm1.13} and \ref{thm1.14} below.

\begin{theorem}\label{thm1.12} (Kleitman, \cite{k}).
Let $n \geq s \geq 3$. Suppose that $\mathcal{K} \subseteq 2^{[n]}$ contains no $s$ pairwise disjoint members. Then the following (i) and (ii) hold:

(i) If $n = sk - 1$, then
\[|\mathcal{K}| \leq \sum_{i \geq k}{{n} \choose {i}};\]

(ii) If  $n = sk + r$ with $0 \leq r \leq s - 2$, then
\[|\mathcal{K}| \leq \sum_{i > k}{{n} \choose {i}} + \frac{s-r-1}{s}{{n} \choose {k}}.\]

\end{theorem}

\hspace{3mm}

We say two subspaces $V_{1}, V_{2}$ in the vector space $\mathbb{F}_{q}^{n}$ are {\em disjoint} if $dim(V_{1} \cap V_{2}) = 0$. The families $\mathcal{V}_{1}, \mathcal{V}_{2}, \dots, \mathcal{V}_{s}$ of subspaces of $\mathbb{F}_{q}^{n}$ are called {\em cross-dependent} if there is no choice of $V_{1} \in \mathcal{V}_{1}, V_{2} \in \mathcal{V}_{2}, \dots, V_{s} \in \mathcal{V}_{s}$ such that $V_{1}, V_{2}, \dots, V_{s}$ are pairwise disjoint.

\begin{theorem}\label{thm1.13}
Let $n \geq s \geq 3$. Suppose that $\mathcal{V} \subseteq \mathcal{L}_{n}(q)$ contains no $s$ pairwise disjoint members. Then the following (i) and (ii) hold:

(i) If $n = sk - 1$, then
\begin{equation}
|\mathcal{V}| \leq \sum_{i \geq k}\qbinom{n}{i};
\end{equation}

(ii)  If $n = sk + r$ with $0 \leq r \leq s - 2$, then
\begin{equation}
|\mathcal{V}| \leq \sum_{i > k}\qbinom{n}{i} + \frac{s - r - 1}{s}\qbinom{n}{k}.
\end{equation}
\end{theorem}
Moreover, the bound in (5) is sharp for the case $n = sk - 1$.

\vspace{3mm}

\begin{theorem}\label{thm1.14}
Suppose that $n = sl + r$ where $0 \leq r < s$ and $s \geq 3$.
Let $\mathcal{V}_{1}, \mathcal{V}_{2}, \dots, \mathcal{V}_{s} \subseteq \mathcal{L}_{n}(q)$ be cross-dependent. Then
\begin{equation}
\sum_{1 \leq i \leq s}|\mathcal{V}_{i}| \leq s \sum_{l < j \leq n}\qbinom{n}{j} + (s-r-1)\qbinom{n}{l}.
\end{equation}
Moreover, the bound in (7) is sharp.
\end{theorem}

\vspace{3mm}

Clearly, Theorem \ref{thm1.14} gives rise to Theorem \ref{thm1.13} by taking $\mathcal{V}_{1} = \mathcal{V}_{2} = \cdots = \mathcal{V}_{s} = \mathcal{V}$.
We will apply Theorem \ref{thm1.9} to derive Theorems \ref{thm1.13} and \ref{thm1.14} in section 4.

\section{Proofs of Theorems 1.9 and 1.11}
First, we introduce some notions and the following covering lemma by Gerbner \cite{g}.
Let $S$ be a set and $S = S_{0} \cup S_{1} \cup \cdots \cup S_{n}$ be a partition (In the discussions here, $S$ is either the set $2^{[n]}$ of all subsets of $[n]$ or the set of all subspaces of $\mathbb{F}_{q}^{n}$, and $S_{i}$ will be level $i$, i.e., $S_{i} = {{[n]} \choose {i}}$ or the set of all $i$-dimensional subspaces of $\mathbb{F}_{q}^{n}$, respectively).
Given a vector $t = (t_{0}, t_{1}, \dots, t_{n})$, we say a family $\Gamma$ of subsets of $S$ is a $t$-{\em covering family} of $S$ if for each $0 \leq i \leq n$, each member in $S_{i}$ is contained in exactly $t_{i}$ sets in the family $\Gamma$.

Given a family $\mathcal{F} \subseteq S$, let $f_{i} = |\mathcal{F} \cap S_{i}|$ and $(f_{0}, f_{1}, \dots, f_{n})$ is called the {\em profile vector} of $\mathcal{F}$.
For a weight $\overline{w} = (w_{0}, w_{1}, \dots, w_{n})$ and a family $\mathcal{F} \subseteq S$,
let $\overline{w}(\mathcal{F}) = \sum_{i = 0}^{n}w_{i}|\mathcal{F} \cap S_{i}|$.
Denote $\overline{w/t} = (w_{0}/t_{0}, w_{1}/t_{1}, \dots, w_{n}/t_{n})$. The following lemma is Lemma 2.1 in \cite{g}.

We say that a property $P$ of families is {\em hereditary} if for any family $\mathcal{F}$ with property $P$,
every subfamily of $\mathcal{F}$ has property $P$. Clearly, forbidding a poset or a configuration is a hereditary property.

\begin{lemma}\label{lem2.1} (Gerbner \cite{g}).
Let $P$ be a hereditary property of subsets (families) of $S$ and $\Gamma$ be a
$t$-covering family of $S$. Assume that there exists a real number $x$ such that for every $G \in \Gamma$, every subset $G'$ of $G$ with property $P$
has $\overline{w/t}(G') \leq x$. Then $\overline{w}(F) \leq |\Gamma|x$ for every $F \subseteq S$ with property $P$.
\end{lemma}

Let $S = \mathcal{L}_{n}(q)$ be the set of all subspaces of $\mathbb{F}_{q}^{n}$.
Denote $S = S_{0} \cup S_{1} \cup \cdots \cup S_{n}$ for a partition of $S$ such that $S_{i}$ is the set of all $i$-dimensional subspaces of $\mathbb{F}_{q}^{n}$.
We construct a $t$-covering family $\Gamma$ of $S$ such that every $\mathcal{G} \in \Gamma$ is a subfamily in $\mathcal{L}_{n}(q)$ isomorphic to Boolean lattice $\mathcal{B}_{n}$ as follows: Choose an arbitrary basis $B = \{v_{1}, \dots, v_{n}\}$ of $\mathbb{F}_{q}^{n}$
and let $\mathcal{G}_{B} = \{span(U) \mid U \subseteq B\}$, i.e., the family of all subspaces that
are generated by subsets of the vectors in $B$.
Obviously the function that maps $H \subseteq [n]$ to the
subspace $span\{v_{x} \mid x \in H\}$ keeps inclusion and intersection properties.
Let $\Gamma$ be the collection of the families $\mathcal{G}_{B}$ over all bases $B$, i.e.,
$\Gamma = \{\mathcal{G}_{B} \mid B \mbox{ is a basis of } \mathbb{F}_{q}^{n}\}$.
Denote
\[\alpha(q, n) = \frac{(q^{n} - 1)(q^{n} - q)(q^{n} - q^{2}) \cdots (q^{n} - q^{n-1})}{n!}.\]
Let $t = (t_{0}, t_{1}, \dots, t_{n})$ be such that for $0 \leq i \leq n$,
\begin{equation}
t_{i} =  \frac{(q^{i} - 1)(q^{i} - q) \cdots (q^{i} - q^{i-1})(q^{n} - q^{i}) \cdots (q^{n} - q^{n-1})}{i!(n-i)!}.
\end{equation}

The next lemma follows from an easy counting.

\begin{lemma}\label{lem2.2}
Let $S = \mathcal{L}_{n}(q)$ be the set of all subspaces of $\mathbb{F}_{q}^{n}$
and $\Gamma = \{\mathcal{G}_{B} \mid B \mbox{ is a basis of } \mathbb{F}_{q}^{n}\}$. Then $|\Gamma| = \alpha(q, n)$ and $\Gamma$ is a $t$-covering of $S$ with
$t = (t_{0}, t_{1}, \dots, t_{n})$ given in (8).
\end{lemma}

\noindent {\bf Proof.} To show $|\Gamma| = \alpha(q, n)$, it suffices to show that $\mathbb{F}_{q}^{n}$ has $\alpha(q, n)$ different bases.
To obtain a basis $B = \{v_{1}, v_{2}, \dots, v_{n}\}$ for $\mathbb{F}_{q}^{n}$, one can choose $n$ nonzero linearly independent vectors in the order
$v_{1}, v_{2}, \dots, v_{n}$ so that $v_{j} \not \in span\{v_{1}, v_{2}, \dots, v_{j-1}\}$ as follows: There are $q^{n} - 1$ ways to choose $v_{1}$, and then there are $q^{n} - q$ ways to choose $v_{2}$, $q^{n} - q^{2}$ ways to choose $v_{3}$, \dots,
$q^{n} - q^{n-1}$ ways to choose $v_{n}$. Thus, there are $(q^{n} - 1)(q^{n} - q)(q^{n} - q^{2}) \cdots (q^{n} - q^{n-1})$ combined ways to form  $B = \{v_{1}, v_{2}, \dots, v_{n}\}$ (in order). Clearly, there are $n!$ ways (permutations on $B$) to form the same $B$ (without order). It follows that there are
\[\alpha(q, n) = \frac{(q^{n} - 1)(q^{n} - q)(q^{n} - q^{2}) \cdots (q^{n} - q^{n-1})}{n!}\]
different bases in $\mathbb{F}_{q}^{n}$.

For any $i$-dimensional subspace $V$ in $\mathbb{F}_{q}^{n}$, similar to the argument above, there are
\[\frac{(q^{i}-1)(q^{i}-q)\cdots (q^{i}-q^{i-1})}{i!}\]
different bases in $V$. For each basis $D$ in $V$, one can extend $D$ to a basis of $\mathbb{F}_{q}^{n}$ in
\[\frac{(q^{n}-q^{i})(q^{n}-q^{i+1})\cdots (q^{n}-q^{n-1})}{(n-i)!}\]
different ways. Therefore, $V$ is contained in
\[t_{i} =  \frac{(q^{i} - 1)(q^{i} - q) \cdots (q^{i} - q^{i-1})(q^{n} - q^{i}) \cdots (q^{n} - q^{n-1})}{i!(n-i)!}\]
different $\mathcal{G}_{B} \in \Gamma$.
It follows that $\Gamma$ is a $t$-covering of $S = \mathcal{L}_{n}(q)$.
$\hfill \Box$

\vspace{3mm}

We now provide a proof for Theorem \ref{thm1.9}.

\vspace{3mm}

\noindent {\bf Proof of Theorem 1.9.}
Assume that $\mathcal{V}$ is a family in $\mathcal{L}_{n}(q)$ which contains no copy of poset $P$.
Let $\Gamma = \{\mathcal{G}_{B} \mid B \mbox{ is a basis of } \mathbb{F}_{q}^{n}\}$. By Lemma \ref{lem2.2},
$\Gamma$ is a $t$-covering of $S = \mathcal{L}_{n}(q)$.
Let $w = (w_{0}, w_{1}, \dots, w_{n})$ with
$w_{i} = \frac{t_{i}}{{{n} \choose {i}}} $ for each $0 \leq i \leq n$.
Then
\[\overline{w/t} = \bigg(\frac{1}{{{n} \choose {0}}}, \frac{1}{{{n} \choose {1}}}, \dots, \frac{1}{{{n} \choose {n}}}\bigg).\]
By the construction, every $\mathcal{G} \in \Gamma$ is isomorphic to Boolean lattice $\mathcal{B}_{n}$. Let $\mathcal{G}' \subseteq \mathcal{G}$ be any subfamily without containing a copy of $P$ and let
\[|\mathcal{G}'| = a_{0}{{n} \choose {0}} + a_{1}{{n} \choose {1}} + \cdots + a_{n}{{n} \choose {n}},\]
where $a_{i}{{n} \choose {i}} = |\mathcal{G}'_{i}|$ and $\mathcal{G}'_{i}$ is the set of all $i$-dimensional subspaces in $\mathcal{G}'$.
Then
\[l(\mathcal{G}') = a_{0} + a_{1} + \cdots + a_{n}.\]
By the assumption,
\[\overline{w/t}(\mathcal{G}') = a_{0}{{n} \choose {0}} \cdot \frac{1}{{{n} \choose {0}}} + a_{1}{{n} \choose {1}} \cdot \frac{1}{{{n} \choose {1}}} + \cdots +
a_{n}{{n} \choose {n}} \cdot \frac{1}{{{n} \choose {n}}} = a_{0} + a_{1} + \cdots + a_{n} = l(\mathcal{G}') \leq c(P).\]
By Lemma \ref{lem2.1}, we have
\begin{equation}
\overline{w}(\mathcal{V})  \leq |\Gamma| \cdot c(P).
\end{equation}
Let
\[|\mathcal{V}| = b_{0}\qbinom{n}{0} + b_{1}\qbinom{n}{1} + \cdots + b_{n}\qbinom{n}{n}\]
with $b_{i}\qbinom{n}{i}$ being the number of $i$-dimensional subspaces in $\mathcal{V}$.
Then
\[l_{q}(\mathcal{V}) = b_{0} +  b_{1} + \cdots + b_{n}.\]
It follows that
\[\overline{w}(\mathcal{V}) = \sum_{i=0}^{n}b_{i}w_{i}\qbinom{n}{i} \hspace{98mm}\]
\[= \sum_{i=0}^{n}\frac{b_{i}}{{{n} \choose {i}}} \frac{(q^{i} - 1)(q^{i} - q) \cdots (q^{i} - q^{i-1})(q^{n} - q^{i}) \cdots (q^{n} - q^{n-1})}{i!(n-i)!}\qbinom{n}{i}  \]
\[= \sum_{i=0}^{n}\frac{b_{i}}{{{n} \choose {i}}} \frac{(q^{i} - 1) \cdots (q^{i} - q^{i-1})(q^{n} - q^{i}) \cdots (q^{n} - q^{n-1})}{i!(n-i)!} \hspace{19mm} \]
\[\hspace{56mm} \cdot \frac{(q^{n} - 1)(q^{n} - q)\cdots (q^{n} - q^{n-1})}{(q^{i} - 1) \cdots (q^{i} - q^{i-1})(q^{n} - q^{i}) \cdots (q^{n} - q^{n-1})}\]
\[ = \sum_{i=0}^{n}\frac{b_{i}}{{{n} \choose {i}}} \frac{\alpha(n, q)n!}{i!(n-i)!} = \sum_{i=0}^{h}\frac{b_{i}}{{{n} \choose {i}}}{{n} \choose {i}}|\Gamma| \hspace{49mm}\]
\[= \sum_{i=0}^{n}b_{i}|\Gamma| = l_{q}(\mathcal{V})|\Gamma|. \hspace{74mm} \]
Combining with (9), we obtain
\[l_{q}(\mathcal{V})|\Gamma| \leq |\Gamma| \cdot c(P), \mbox{ that is, } l_{q}(\mathcal{V}) \leq c(P).\]
$\hfill \Box$

\vspace{3mm}

\noindent {\bf Proof of Theorem 1.11.} The proof is obtained by modifying the proof of Theorem \ref{thm1.9} above. Indeed, let $w = (w_{0}, w_{1}, \dots, w_{n})$ with $w_{i} = \frac{\beta_{i}t_{i}}{{{n} \choose {i}}}$. We have
\[\overline{w/t} = \bigg(\frac{\beta_{0}}{{{n} \choose {0}}}, \frac{\beta_{1}}{{{n} \choose {1}}}, \dots, \frac{\beta_{n}}{{{n} \choose {n}}}\bigg),\]
\[l^{*}(\mathcal{G}', \beta) = \sum_{i = 0}^{n}\frac{\beta_{i}|\mathcal{G}'_{i}|}{{{n} \choose {i}}} = a_{0}\beta_{0} + a_{1}\beta_{1} + \cdots + a_{n}\beta_{n},\]
\[\overline{w/t}(\mathcal{G}') = a_{0}{{n} \choose {0}} \cdot \frac{\beta_{0}}{{{n} \choose {0}}} + a_{1}{{n} \choose {1}} \cdot \frac{\beta_{1}}{{{n} \choose {1}}} + \cdots + a_{n}{{n} \choose {n}} \cdot \frac{\beta_{n}}{{{n} \choose {n}}} \hspace{8mm}\]
\[= a_{0}\beta_{0} + a_{1}\beta_{1} + \cdots + a_{n}\beta_{n} = l^{*}(\mathcal{G}', \beta) \leq c(P).\]
By Lemma \ref{lem2.1}, we have
\begin{equation}
\overline{w}(\mathcal{V})  \leq |\Gamma| \cdot c(P).
\end{equation}
On the other hand, we have
\[\overline{w}(\mathcal{V}) = \sum_{i=0}^{n}b_{i}w_{i}\qbinom{n}{i} \hspace{100mm}\]
\[= \sum_{i=0}^{n}\frac{b_{i}\beta_{i}}{{{n} \choose {i}}} \frac{(q^{i} - 1)(q^{i} - q) \cdots (q^{i} - q^{i-1})(q^{n} - q^{i}) \cdots (q^{n} - q^{n-1})}{i!(n-i)!}\qbinom{n}{i}  \]
\[= \sum_{i=0}^{n}b_{i}\beta_{i}|\Gamma| = l^{*}_{q}(\mathcal{V}, \beta)|\Gamma|. \hspace{68mm} \]
Together with (10), it follows that
\[l^{*}_{q}(\mathcal{V}, \beta) \leq c(P).\]
$\hfill \Box$

\section{Sperner families and Diamond-free families}
In this section, we apply Theorems \ref{thm1.9} and \ref{thm1.11} to derive generalizations of results on Sperner families and diamond-free families of subsets of $[n]$.

\subsection{Sperner and $k$-Sperner families}
Clearly, Theorems \ref{thm1.5} and \ref{thm1.9} give the next generalization of Theorem \ref{thm1.5} which implies Theorem \ref{thm1.3}.

\begin{theorem}\label{thm3.1} (q-Analogue LYM inequality).  Assume that $\mathcal{V}$ is a Sperner family
of subspaces of $\mathbb{F}_{q}^{n}$. Then
\[l_{q}(\mathcal{V}) = \sum_{V \in \mathcal{V}}\frac{1}{\qbinom{n}{|V|}} \leq 1.\]
\end{theorem}

Also, Theorems \ref{thm1.6} and \ref{thm1.9} give the following generalization of Theorem \ref{thm1.6} which implies Theorem \ref{thm1.4}.

\begin{theorem}\label{thm3.2}  Assume that $\mathcal{V}$ is a $k$-Sperner family
of subspaces of $\mathbb{F}_{q}^{n}$. Then
\[l_{q}(\mathcal{V}) = \sum_{V \in \mathcal{V}}\frac{1}{\qbinom{n}{|V|}} \leq k.\]
\end{theorem}

\subsection{Diamond-free families}
The diamond poset, denoted by $\mathcal{Q}_{2}$, is defined on four elements
$x, y, z, w$ with the relations $x < y, z$ and $y, z < w$. Recall that
$ex(n, P)$ is the maximum size of a family of subsets of
$[n]$ not containing $P$, $ex(n, \mathcal{Q}_{2})$ is the maximum size of a diamond-free family of subsets of
$[n]$. Despite of efforts made by many researchers, even the asymptotic value of $ex(n, \mathcal{Q}_{2})$ has yet
to be determined. It is conjectured that
\[ex(n, \mathcal{Q}_{2}) = (2 + o(1)){{n} \choose {\lfloor \frac{n}{2}\rfloor}}.\]
Note that the two middle levels of the Boolean lattice do not contain a diamond, we have
\[ex(n, \mathcal{Q}_{2}) \geq (2 - o(1)){{n} \choose {\lfloor \frac{n}{2}\rfloor}}.\]
Czabarka et al. \cite{cdjs} gave infinitely many asymptotically tight
constructions by using random set families defined from posets based on Abelian groups.
Using an elegant argument, Griggs et al. \cite{gll} showed that $ex(n, \mathcal{Q}_{2}) < 2.296{{n} \choose {\lfloor \frac{n}{2}\rfloor}}$.
This bound was further improved to $(2.25 + o(1)){{n} \choose {\lfloor \frac{n}{2}\rfloor}}$ by Kramer et al. \cite{kmy}.
The best known upper bound on $ex(n, \mathcal{Q}_{2})$ is $(2.20711 + o(1)){{n} \choose {\lfloor \frac{n}{2}\rfloor}}$
provided by Gr\'osz et al. \cite{gmt} as follows.

\begin{theorem}\label{thm3.3} (Gr\'osz, Methuku, and Tompkins, \cite{gmt}).
\[ex(n, \mathcal{Q}_{2}) \leq \bigg(\frac{\sqrt{2}+3}{2} + o(1)\bigg){{n} \choose {\lfloor \frac{n}{2}\rfloor}} < (2.20711 + o(1)){{n} \choose {\lfloor \frac{n}{2}\rfloor}}.\]
\end{theorem}

In the proof of Theorem \ref{thm3.3} which is Theorem 1.15 in \cite{gmt}, Gr\'osz et al. derived the next upper bound for the Lubell function $l(\mathcal{F})$.

\begin{proposition}\label{prop3.4} (Gr\'osz, Methuku, and Tompkins, \cite{gmt}).
Suppose that $\mathcal{F} \subseteq 2^{[n]}$ is a family which contains no $\mathcal{Q}_{2}$. Then
\[l(\mathcal{F}) \leq \frac{\sqrt{2}+3}{2} + o(1) < 2.20711 + o(1).\]
\end{proposition}

Clearly, Theorem \ref{thm1.9} and Proposition \ref{prop3.4} imply immediately the following generalization of Proposition \ref{prop3.4}.

\begin{proposition}\label{prop3.5}
Suppose that $\mathcal{V}$ is a family of subspaces of $\mathbb{F}_{q}^{n}$ which contains no $\mathcal{Q}_{2}$. Then
\[l_{q}(\mathcal{F}) \leq \frac{\sqrt{2}+3}{2} + o(1) < 2.20711 + o(1).\]
\end{proposition}

From Proposition \ref{prop3.5}, the next generalization of Theorem \ref{thm3.3} follows directly, where $ex_{q}(n, \mathcal{Q}_{2})$ denotes the maximum size of a family of subspaces of $\mathbb{F}_{q}^{n}$ which contains no diamond $\mathcal{Q}_{2}$.

\begin{theorem}\label{thm3.6}
\[ex_{q}(n, \mathcal{Q}_{2}) \leq \bigg(\frac{\sqrt{2}+3}{2} + o(1)\bigg)\qbinom{n}{\lfloor \frac{n}{2}\rfloor}
< (2.20711 + o(1)) \qbinom{n}{\lfloor \frac{n}{2}\rfloor}.\]
\end{theorem}

Note that the two middle levels of the lattice $\mathcal{L}_{n}(q)$ of subspaces of $\mathbb{F}_{q}^{n}$ do not contain a diamond $\mathcal{Q}_{2}$.  The maximum size $ex_{q}(n, \mathcal{Q}_{2})$ of a family of subspaces of $\mathbb{F}_{q}^{n}$ containing no $\mathcal{Q}_{2}$ satisfies
\[ex_{q}(n, \mathcal{Q}_{2}) \geq (2 - o(1))\qbinom{n}{\lfloor \frac{n}{2}\rfloor}.\]
We propose the following conjecture analogous to the corresponding one for Boolean lattices.

\begin{conjecture}\label{conj3.7}
The maximum size $ex_{q}(n, \mathcal{Q}_{2})$ of a family of subspaces of $\mathbb{F}_{q}^{n}$ containing no diamond $\mathcal{Q}_{2}$ satisfies
\[ex_{q}(n, \mathcal{Q}_{2}) = 2 \qbinom{n}{\lfloor \frac{n}{2}\rfloor}.\]
\end{conjecture}

\subsection{Sharpening the LYM inequality}
Note that
\[l(\mathcal{F}) = \sum_{F \in \mathcal{F}}\frac{1}{{{n} \choose {|F|}}} = \sum_{i = 0}^{n}\frac{|\mathcal{F}_{i}|}{{{n} \choose {i}}}\]
with each $\mathcal{F}_{i} = \mathcal{F} \cap {{[n]} \choose {i}}$.
The next result in \cite{efkss} sharpens
the LYM inequality by increasing the coefficients  $\frac{1}{{{n} \choose {i}}}$.

\begin{theorem}\label{thm3.8} (Erd\H{o}s, Frankl, Kleitman, Saks, and Sz$\acute{e}$kely, \cite{efkss}).
Let $\mathcal{F} \subseteq 2^{[n]}$ is a Sperner family (or antichain) and let $k$ be the
smallest integer $t$ for which $\sum_{i \leq t}\frac{|\mathcal{F}_{i}|}{{{n-1} \choose {i-1}}} > 1$.
Then
\[ \sum_{i < k}\frac{k}{i}\frac{|\mathcal{F}_{i}|}{{{n} \choose {i}}} + \sum_{i \geq k}\frac{n-k}{n-i}\frac{|\mathcal{F}_{i}|}{{{n} \choose {i}}} \leq 1.\]
\end{theorem}

By Theorems \ref{thm1.11} and \ref{thm3.8}, we have the next generalization of Theorem \ref{thm3.8} directly.

\begin{theorem}\label{thm3.9}
Assume that $\mathcal{V}$ is a Sperner family of subspaces of $\mathbb{F}_{q}^{n}$ and let $k$ be the
smallest integer $t$ for which $\sum_{i \leq t}\frac{|\mathcal{V}_{i}|}{\qbinom{n-1}{i-1}} > 1$, where $\mathcal{V}_{i} = \mathcal{V} \cap \qbinom{[n]}{i}$.
Then
\[ \sum_{i < k}\frac{k}{i}\frac{|\mathcal{V}_{i}|}{\qbinom{n}{i}} + \sum_{i \geq k}\frac{n-k}{n-i}\frac{|\mathcal{V}_{i}|}{\qbinom{n}{i}} \leq 1.\]
\end{theorem}

Recently, Malec and Tompkins \cite{mt} proved the following strengthening of Theorem \ref{thm1.6}.
For a family  $\mathcal{F} \subseteq 2^{[n]}$, let
\[c(F) = \max \{k \mid F \mbox{ participates in a } k\mbox{-chain consisting of sets from } \mathcal{F}\}.\]

\begin{theorem}\label{thm3.10} (Malec and Tompkins, \cite{mt})
Let $\mathcal{F} \subseteq 2^{[n]}$ be an arbitrary family of sets, then
\[\sum_{F \in \mathcal{F}}\frac{1}{c(F){{n} \choose {|F|}}} \leq 1.\]
\end{theorem}
Obviously, if a family $\mathcal{F} \subseteq 2^{[n]}$ is $k$-Sperner, then $c(F) \leq k$ for every $F \in \mathcal{F}$, and so Theorem \ref{thm3.10}
implies Theorem \ref{thm1.6}.
By applying Theorems \ref{thm1.11} and \ref{thm3.10}, we obtain the following generalization of Theorem \ref{thm3.10},
where
\[c_{q}(V) = \max \{k \mid V \mbox{ participates in a } k\mbox{-chain consisting of subspaces from } \mathcal{V}\}.\]

\begin{theorem}\label{thm3.11}
Let $\mathcal{V}$ be a family of subspaces of $\mathbb{F}_{q}^{n}$, then
\[\sum_{V \in \mathcal{V}}\frac{1}{c_{q}(V)\qbinom{n}{dim(V)}} \leq 1.\]
\end{theorem}

\section{Families with no s pairwise disjoint members}
In this section, we will apply Theorem \ref{thm1.9} to give proofs for Theorems \ref{thm1.13} and \ref{thm1.14}.
For non-triviality, we assume that $q \geq 2$.

\subsection{Preliminaries}
Given a family $\mathcal{F} \subseteq 2^{[n]}$, Frankl \cite{f} calls the Lubell function $l(\mathcal{F})$ the {\em binomial norm}, and denoted by
\[||\mathcal{F}|| = \sum_{F \in \mathcal{F}}\frac{1}{{{n} \choose {|F|}}} = l(\mathcal{F}).\]
For $0 \leq i \leq n$, set
\[\mathcal{F}_{i} = \{ F \in \mathcal{F} \mid |F| = i\}\]
and define
\[\varphi(i) = \frac{|\mathcal{F}_{i}|}{{{n} \choose {i}}}.\]
Then
\[||\mathcal{F}|| = \sum_{i = 0}^{n}\frac{|\mathcal{F}_{i}|}{{{n} \choose {i}}} = \sum_{i = 0}^{n}\varphi(i).\]
Define
\[\varrho_{n}(\mathcal{F}) = \frac{||\mathcal{F}||}{n+1}.\]
Clearly, the binomial norm satisfies
\[0 \leq ||\mathcal{F}|| \leq n + 1\]
which implies that
\[0 \leq \varrho_{n}(\mathcal{F}) \leq 1.\]

\begin{definition}\label{def4.1}
The families $\mathcal{F}_{1}, \mathcal{F}_{2}, \dots, \mathcal{F}_{s}$ in $2^{[n]}$ are called {\em cross-dependent} if there is no choice of
$F_{1} \in \mathcal{F}_{1}, F_{2} \in \mathcal{F}_{2}, \dots, F_{s} \in \mathcal{F}_{s}$ such that $F_{1}, F_{2}, \dots, F_{s}$ are pairwise disjoint.
\end{definition}

The following theorem is Theorem 1.6 in \cite{f}.

\begin{theorem}\label{thm4.2} (Frankl, \cite{f}).
Suppose that $\mathcal{F}_{1}, \mathcal{F}_{2}, \dots, \mathcal{F}_{s} \subseteq 2^{[n]}$ are cross-dependent. Then
\begin{equation}
\varrho_{n}(\mathcal{F}_{1}) + \varrho_{n}(\mathcal{F}_{2}) + \cdots + \varrho_{n}(\mathcal{F}_{s}) \leq s - 1.
\end{equation}
\end{theorem}

Obviously, by the definition, (11) is equivalent to
\begin{equation}
||\mathcal{F}_{1}|| + ||\mathcal{F}_{2}|| + \cdots + ||\mathcal{F}_{s}|| \leq (s - 1)(n + 1).
\end{equation}

By taking $\mathcal{F}_{1} = \mathcal{F}_{2} = \dots = \mathcal{F}_{s} = \mathcal{F}$, Theorem \ref{thm4.2} has the next useful corollary.

\begin{corollary}\label{cor4.3} (Frankl, \cite{f}).
Suppose that $\mathcal{F} \subseteq 2^{[n]}$ contains no $s$ pairwise disjoint members. Then
\[||\mathcal{F}|| \leq \frac{(s-1)(n+1)}{s}. \]
\end{corollary}

Recall that two subspaces $V_{1}, V_{2}$ in the vector space $\mathbb{F}_{q}^{n}$ are disjoint if $dim(V_{1} \cap V_{2}) = 0$.

\begin{definition}\label{def4.4}
The families $\mathcal{V}_{1}, \mathcal{V}_{2}, \dots, \mathcal{V}_{s}$ of subspaces of $\mathbb{F}_{q}^{n}$ are called {\em cross-dependent} if there is no choice of
$V_{1} \in \mathcal{V}_{1}, V_{2} \in \mathcal{V}_{2}, \dots, V_{s} \in \mathcal{V}_{s}$ such that $V_{1}, V_{2}, \dots, V_{s}$ are pairwise disjoint.
\end{definition}

Let $\mathcal{V}$ be a family of subspaces of the vector space $\mathbb{F}_{q}^{n}$. We define the $q$-binomial norm of $\mathcal{V}$ to be
\[||\mathcal{V}||_{q} = l_{q}(\mathcal{V}) = \sum_{V \in \mathcal{V}}\frac{1}{\qbinom{n}{|V|}}.\]
For $0 \leq i \leq n$, set
\[\mathcal{V}_{i} = \{ V \in \mathcal{V} \mid dim(V) = i\}\]
and define
\[\varphi_{q}(i) = \frac{|\mathcal{V}_{i}|}{\qbinom{n}{i}}.\]
Then
\begin{equation}
||\mathcal{V}||_{q} = \sum_{i = 0}^{n}\frac{|\mathcal{V}_{i}|}{\qbinom{n}{i}} = \sum_{i = 0}^{n}\varphi_{q}(i).
\end{equation}
By Theorems \ref{thm1.9} and \ref{thm4.2}, Corollary \ref{cor4.3}, and (12), we obtain the following generalizations of Theorem \ref{thm4.2} and Corollary \ref{cor4.3} immediately.

\begin{theorem}\label{thm4.5}
Suppose that families $\mathcal{V}_{1}, \mathcal{V}_{2}, \dots, \mathcal{V}_{s}$ of subspaces of $\mathbb{F}_{q}^{n}$ are cross-dependent.
Then
\[||\mathcal{V}_{1}||_{q} + ||\mathcal{V}_{2}||_{q} + \cdots + ||\mathcal{V}_{s}||_{q} \leq (s - 1)(n+1).\]
\end{theorem}

\begin{theorem}\label{thm4.6}
Suppose that $\mathcal{V}$ is a family of subspaces of $\mathbb{F}_{q}^{n}$ which contains no $s$ pairwise disjoint members. Then
\begin{equation}
||\mathcal{V}||_{q} \leq \frac{(s-1)(n+1)}{s}.
\end{equation}
\end{theorem}

\subsection{Proofs of Theorems 1.13 and 1.14}
Next, we derive Theorems \ref{thm1.13} and \ref{thm1.14} by applying Theorems \ref{thm4.5} and \ref{thm4.6},
following an approach similar to that of Frankl \cite{f}, but the analysis is more complicated. Theorem \ref{thm1.9}
plays a crucial role in deriving Theorems \ref{thm4.5} and \ref{thm4.6}.

The shadow of a family of subsets is extended naturally to vector spaces as follows: For a family $\mathcal{V} \subseteq \qbinom{[n]}{k}$, the shadow of  $\mathcal{V}$
is defined to be
\[\partial\mathcal{V} = \bigg\{ W \in \qbinom{[n]}{k-1} \bigg| W \subseteq V \in \mathcal{V}\bigg\}.\]
The next $q$-analogue of Lov$\acute{a}$sz's theorem for the Boolean lattices is Theorem 1.4 in \cite{cp}.

\begin{theorem}\label{thm4.7} (Chowdhury and  Patk\'os, \cite{cp})
Let $\mathcal{V} \subseteq \qbinom{[n]}{k}$ and let $y$ be the real number satisfying $|\mathcal{V}| = \qbinom{y}{k}$.
Then
\[|\partial\mathcal{V}| \geq \qbinom{y}{k-1}.\]
\end{theorem}

A family $\mathcal{V}$ of subspaces of $\mathbb{F}_{q}^{n}$ is said to be a {\em complex} if $W \subseteq V \in \mathcal{V}$ implies $W \in \mathcal{V}$.
We have the next lemma.

\begin{lemma}\label{lem4.8} Let a family $\mathcal{V}$ of subspaces of $\mathbb{F}_{q}^{n}$ be a complex and let
$\mathcal{V}_{i} = \{ V \in \mathcal{V} \mid dim(V) = i\}$ for $0 \leq i \leq n$.
Then the function $\varphi_{q}(i) = \frac{|\mathcal{V}_{i}|}{\qbinom{n}{i}}$ is monotonic decreasing, i.e.,
\[\varphi_{q}(0) \geq \varphi_{q}(1) \geq \cdots \geq \varphi_{q}(n).\]
\end{lemma}

\vspace{3mm}

\noindent {\bf Proof.}
Note that
\[\qbinom{y}{k-1} = \frac{q^{k} - 1}{q^{y - k + 1} - 1}\qbinom{y}{k},\]
\[\qbinom{n}{k-1} = \frac{q^{k} - 1}{q^{n - k + 1} - 1}\qbinom{n}{k}.\]
It follows that for any $y \leq n$,
\[\frac{\qbinom{y}{k}}{\qbinom{n}{k}} \leq \frac{\qbinom{y}{k-1}}{\qbinom{n}{k-1}}.\]
Let $1 \leq k \leq n$ be fixed and let $y$ be the real number satisfying $|\mathcal{V}_{k}| = \qbinom{y}{k}$. Then $y \leq n$. By Theorem \ref{thm4.7},
\[|\partial\mathcal{V}_{k}| \geq \qbinom{y}{k-1}.\]
Since $\mathcal{V}$ is a complex, we have $\partial\mathcal{V}_{k} \subseteq \mathcal{V}_{k-1}$.
It follows that
\[\varphi_{q}(k-1) = \frac{|\mathcal{V}_{k-1}|}{\qbinom{n}{k-1}} \geq \frac{|\partial\mathcal{V}_{k}|}{\qbinom{n}{k-1}} \geq \frac{\qbinom{y}{k-1}}{\qbinom{n}{k-1}}
\geq \frac{\qbinom{y}{k}}{\qbinom{n}{k}} = \frac{|\mathcal{V}_{k}|}{\qbinom{n}{k}} = \varphi_{q}(k).\]
$\hfill \Box$

\vspace{3mm}

Note that $\qbinom{n}{k}$ defined by (1) satisfies the symmetric and unimodal properties similar to
those of the binomial coefficients ${{n} \choose {k}}$. We have
the following fact.

\begin{lemma}\label{lem4.9}
Let $n, k, l, q$ be positive integers such that $l < k \leq n -1$, $q \geq 2$, and $n \geq 3l$. Then
\begin{equation}
\sum_{l \leq j \leq k}\qbinom{n}{j} \geq (k - l + 1)\qbinom{n}{l}.
\end{equation}
\end{lemma}

\noindent {\bf Proof.}
By the symmetric and unimodal properties of $\qbinom{n}{i}$, we have $\qbinom{n}{l} < \qbinom{n}{j}$ for $l < j < n - l$. Thus, (15) is obviously true for $k \leq n - l$, and so the lemma holds for $l = 0, 1$. Hence we assume $l \geq 2$ and $k > n - l$.

Observe that
\begin{equation}
\qbinom{n}{l+1} = \frac{q^{n-l} - 1}{q^{l+1} - 1}\qbinom{n}{l},
\end{equation}
\begin{equation}
\qbinom{n}{l+2} = \frac{q^{n-l-1} - 1}{q^{l+2} - 1}\qbinom{n}{l+1}.
\end{equation}
We claim that $\qbinom{n}{l+1} \geq 2\qbinom{n}{l}$ and $\qbinom{n}{l+2} \geq 2\qbinom{n}{l}$ when $q \geq 2$ and $n \geq 3l$ except for $n = 3l$ and $l = 2$.
To see the claim, by (16) and (17), it suffices to show that (i) $q^{n-l} - 1 \geq 2(q^{l+1} - 1)$ and (ii) $q^{n - l - 1} - 1 \geq q^{l+2} - 1$.
For (ii), it is equivalent to $n - l - 1 \geq l + 2$, i.e., $n \geq 2l +3$ which is true if $n \geq 3l$ and $l \geq 3$ or $n \geq 3l + 1$ and $l = 2$.
For (i), since $q \geq 2$ and $l \geq 2$,  $n \geq 3l$ implies $n \geq 2l + 2$ and so $q^{n-l} \geq q^{l+2}$
which implies that $q^{n-l} - 1 \geq q^{l+2} - 1 \geq 2q^{l+1} - 1 > 2(2^{l+1} - 1)$. Thus, the claim holds.

We consider the following two cases:

\noindent {\bf Case 1.} $n \neq 6$ when $l = 2$.
In this case, by the claim, we have
\begin{equation}
\qbinom{n}{l+1} + \qbinom{n}{l+2} \geq 4\qbinom{n}{l}.
\end{equation}
Since the terms $\qbinom{n}{j}$ that exceed $\qbinom{n}{l}$ are $\qbinom{n}{j}$ for $l+1 \leq j \leq n-l-1$ (there are $n - 2l - 1$ such terms) and $\qbinom{n}{l+1}$ and $\qbinom{n}{l+2}$ are the smallest, it follows (18) that
\begin{equation}
\sum_{l+1 \leq j \leq n-l-1}\qbinom{n}{j} \geq 2(n - 2l - 1)\qbinom{n}{l}.
\end{equation}
Adding $\qbinom{n}{l}$ and $\qbinom{n}{n-l}$ to (19) gives
\begin{equation}
\sum_{l \leq j \leq n-l}\qbinom{n}{j} \geq 2(n - 2l)\qbinom{n}{l}.
\end{equation}
Since $n \geq 3l$ and $k \leq n - 1$, $2(n - 2l) \geq k - l + 1$. It follows from (20) that
\begin{equation}
\sum_{l \leq j \leq n-l}\qbinom{n}{j} \geq (k - l + 1)\qbinom{n}{l}.
\end{equation}
Clearly, (21) implies (15) when $k > n - l$.

\noindent {\bf Case 2.} $n = 6$ and $l = 2$. Since $n - l < k \leq n - 1$, we have $k = n -1 = 5$. To prove (15), we need to verify
\begin{equation}
\sum_{2 \leq j \leq 5}\qbinom{6}{j} \geq 4\qbinom{6}{2}.
\end{equation}
By symmetry, $\qbinom{6}{2} = \qbinom{6}{4}$. Hence (22) becomes
\begin{equation}
\qbinom{6}{3} + \qbinom{6}{5} \geq 2\qbinom{6}{2}.
\end{equation}
Since $q \geq 2$, $q^{3} \geq 2q^{2} \geq q^{2} + 2q > q^{2} + q + 1$.
It follows from (16) that
\[\qbinom{6}{3} = \frac{q^{4} - 1}{q^{3} - 1}\qbinom{6}{2} = \frac{q^{3} + q^{2} + q + 1}{q^{2} + q + 1} \geq 2\qbinom{6}{2}\]
which implies (23).
Therefore, the lemma holds.
$\hfill \Box$

\hspace{3mm}

Recall that $\qbinom{n}{k}$ is equal to the number of all $k$-dimensional
subspaces of $\mathbb{F}_{q}^{n}$. We have
\begin{equation}
|\mathcal{L}_{n}(q)| = \sum_{i=0}^{n}\qbinom{n}{i}
\end{equation}
which implies $||\mathcal{V}||_{q} \leq n + 1$ for any $\mathcal{V} \subseteq \mathcal{L}_{n}(q)$ (see (13)).

\begin{definition}\label{def4.10}
Fix a pair $(n, l)$ with $0 \leq l < n$. We say $(n, l)$ is $q$-{\em perfect} if the following inequality holds for all families $\mathcal{V}$ of subspaces of $\mathbb{F}_{q}^{n}$ which are complexes with $||\mathcal{V}||_{q} < n+1$:
\begin{equation}
|\mathcal{V}| \geq \sum_{0 \leq i < l}\qbinom{n}{i} + (||\mathcal{V}||_{q} - l)\qbinom{n}{l}.
\end{equation}
\end{definition}

The next proposition is a $q$-analogue of Proposition 4.2 in \cite{f}, with a similar proof by applying Lemmas \ref{lem4.8} and \ref{lem4.9} and replacing
${{n} \choose {i}}$'s by $\qbinom{n}{i}$'s (see Appendix 1 for a detailed proof).

\begin{proposition}\label{prop4.11}
If $n \geq 3l$, then $(n, l)$ is $q$-perfect.
\end{proposition}

We are now ready to provide the following proofs for Theorems \ref{thm1.13} and \ref{thm1.14} by applying Theorems \ref{thm4.5} and \ref{thm4.6} and Proposition \ref{prop4.11}.

\vspace{3mm}

\noindent {\bf Proof of Theorem 1.13.}
Note that if $\mathcal{V}$ contains no s pairwise disjoint members, then the up-set
\[\mathcal{V}^{*} = \{W \in \mathcal{L}_{n}(q) \mid \exists V \in \mathcal{V}, \mbox{ } V \subseteq W\}\]
has the same property as well. Thus, we may assume that $\mathcal{V}$ itself is an up-set.
Define
\[\mathcal{F}  = \mathcal{L}_{n}(q) \setminus \mathcal{V}.\]
Then $\mathcal{F}$ is a complex as $\mathcal{V}$ is an up-set. By Theorem \ref{thm4.6}, we have
\[||\mathcal{V}||_{q} \leq \frac{(s-1)(n+1)}{s}\]
which implies that $||\mathcal{F}||_{q} \geq \frac{n+1}{s}$ as $ ||\mathcal{L}_{n}(q)||_{q} = n+1$ by (13) and (24).

For (i), we have $n = sk - 1$ and so $||\mathcal{F}||_{q} \geq k$.
Apply Proposition \ref{prop4.11} with $l = k - 1$, we obtain
\[|\mathcal{F}| \geq \sum_{0 \leq i < l}\qbinom{n}{i} + (||\mathcal{F}||_{q} - l)\qbinom{n}{l} \geq \sum_{0 \leq i \leq k-1}\qbinom{n}{i}\]
which gives (5).

For (ii), we have  $n = sk + r$ with $0 \leq r \leq s - 2$ and so $||\mathcal{F}||_{q} \geq k + \frac{r + 1}{s}$.
Apply Proposition \ref{prop4.11} with $l = k$, we obtain
\[|\mathcal{F}| \geq \sum_{0 \leq i < k}\qbinom{n}{i} + (||\mathcal{F}||_{q} - k)\qbinom{n}{k} \geq \sum_{0 \leq i \leq k-1}\qbinom{n}{i} + \frac{r+1}{s}\qbinom{n}{k}\]
which gives (6).

For the sharpness in the case $n = sk - 1$, take the family $\mathcal{V}$ of all subspaces of dimension greater than or equal to $k$. Then
$|\mathcal{V}| = \sum_{i \geq k}\qbinom{n}{i}$ and $\mathcal{V}$ contains no s pairwise disjoint members.

\hspace{60mm} $\hfill \Box$

\vspace{3mm}

\noindent {\bf Proof of Theorem 1.14.}
Set $\mathcal{W}_{i} = \mathcal{L}_{n}(q) \setminus \mathcal{V}_{i}$ for $1 \leq i \leq s$. By (24) and $n = sl + r$, (7) is equivalent to
\begin{equation}
\sum_{1 \leq i \leq s}|\mathcal{W}_{i}| \geq s \sum_{0 \leq j < l}\qbinom{n}{j} + (r+1)\qbinom{n}{l}.
\end{equation}
By Theorem \ref{thm4.5}, we obtain
\[\sum_{1 \leq i \leq s}||\mathcal{V}_{i}||_{q} \leq (s-1)(n+1).\]
Note that $||\mathcal{W}_{i}||_{q} = (n + 1) - ||\mathcal{V}_{i}||_{q}$ for each $1 \leq i \leq s$. It follows that
\[\sum_{1 \leq i \leq s}||\mathcal{W}_{i}||_{q} \geq n + 1.\]
Since $(n, l)$ is $q$-perfect by Proposition \ref{prop4.11}, (25) implies that
\[\sum_{1 \leq i \leq s}|\mathcal{W}_{i}| \geq s \sum_{0 \leq j < l}\qbinom{n}{j} + (n+1 - sl)\qbinom{n}{l},\]
which gives (26).

For the sharpness of the bound in (7), consider the families
\[\mathcal{V}_{1} = \cdots = \mathcal{V}_{r+1} = \{V \in \mathcal{L}_{n}(q) \mid dim(V) > l\},\]
\[\mathcal{V}_{r+2} = \cdots = \mathcal{V}_{s} = \{V \in \mathcal{L}_{n}(q) \mid dim(V) \geq l\}.\]
This shows that the bound in (7) is best possible.
$\hfill \Box$

\section{Families without $d$-dimensional $q$-algebras in linear lattices}
In this section, we apply Theorem  \ref{thm1.9} to derive generalizations of interesting results about the maximum size of a family $\mathcal{F} \subseteq  2^{[n]}$ which does not contain a $d$-dimensional Boolean algebra.
Recall that $2^{[n]}$ is the power set of $[n] = \{1, 2, \dots, n\}$. The following concept of $d$-dimensional Boolean algebra is given in \cite{grs} and \cite{jlm} (also called $d$-dimensional combinatorial subspace in \cite{p}).

\begin{definition}\label{def5.1}
A collection $\mathcal{B} \subseteq  2^{[n]}$ forms a $d$-{\em dimensional Boolean algebra} if there exist
pairwise disjoint sets $X_{0}, X_{1}, \dots, X_{d} \subseteq [n]$, all non-empty with perhaps the exception of $X_{0}$,
so that
\[\mathcal{B} = \bigg\{X_{0} \cup(\cup_{i \in I}X_{i}) \mid I \subseteq [d]\bigg\}.\]
\end{definition}

\begin{definition}\label{def5.2}
Given a positive integer d,
define $b(n, d)$ to be the maximum size of a family $\mathcal{F} \subseteq  2^{[n]}$ which does
not contain a $d$-dimensional Boolean algebra ($\mathcal{B}_{d}$-free).
\end{definition}

Note that a $1$-dimensional Boolean algebra is simply a pair of sets, one
contained in the other and so, by Sperner's theorem (Theorem \ref{thm1.1}),
\[b(n, 1) = {{n} \choose {\lfloor \frac{n}{2}\rfloor}} \sim \sqrt{\frac{2}{\pi}}n^{-\frac{1}{2}}\cdot 2^{n} = o(2^{n}).\]
Erd\H{o}s and Kleitman \cite{ek} found that there exist constants $c_{1}$ and $c_{2}$ so that
for $n$ sufficiently large,
\[c_{1}n^{-\frac{1}{4}}\cdot 2^{n} \leq b(n, 2) \leq c_{2}n^{-\frac{1}{4}}\cdot 2^{n}.\]
In \cite{grs}, the following bounds on $b(n, d)$ are proved by Gunderson et al.:
\begin{equation}
n^{-\frac{(1 + o(1))d}{2^{d+1}-2}}\cdot 2^{n} \leq b(n, d) \leq 10^{d}2^{-2^{1-d}}d^{d-2^{-d}}n^{-\frac{1}{2^{d}}}\cdot 2^{n} = o(2^{n}).
\end{equation}
The upper bound in (27) is improved by Johnston et al. \cite{jlm} and Polymath \cite{p} as follows.

\begin{theorem}\label{thm5.3} (Polymath, \cite{p}).
The maximum size $b(n, d)$ of a family $\mathcal{F} \subseteq  2^{[n]}$ which does
not contain a $d$-dimensional Boolean algebra satisfies
\[b(n, d) \leq \bigg(\frac{25}{n}\bigg)^{\frac{1}{2^{d}}}\cdot 2^{n}.\]
\end{theorem}

We remark that Sperner's theorem and Theorem \ref{thm5.3} were applied in a new proof of Furstenberg and Katznelson's
density Hales-Jewett theorem (which implies the well-known Szemer\'edi theorem, see Theorem \ref{thm6.1} below) by Polymath \cite{p}.

\begin{theorem}\label{thm5.4} (Johnston, Lu, and Milans, \cite{jlm}).
There is a positive constant $C$, independent of $d$, such that for every $d$ and all sufficiently large $n$, the following is true.
\[b(n, d) \leq Cn^{-\frac{1}{2^{d}}}\cdot 2^{n}.\]
\end{theorem}

In order to prove Theorem \ref{thm5.4}, Johnston et al. \cite{jlm} obtained the following useful result as Corollary 1, where $l(\mathcal{F})$ is the Lubell function defined in (2).

\begin{theorem}\label{thm5.5} (Johnston, Lu, and Milans, \cite{jlm}).
For $d \geq 3$ and $n \geq (2^{d} - \frac{2}{ln2})^{2}$, every family $\mathcal{F} \subseteq  2^{[n]}$ containing no $d$-dimensional Boolean algebra satisfies
\[l(\mathcal{F}) \leq 2(n+1)^{1 - 2^{1-d}}.\]
\end{theorem}

Similar to Definitions \ref{def5.1} and \ref{def5.2}, one can define the corresponding concept and value for vector spaces as follows.
Recall that $\mathcal{L}_{n}(q)$ denotes the linear lattice of the set of all subspaces of $\mathbb{F}_{q}^{n}$.

\begin{definition}\label{def5.6}
A collection $\mathcal{B} \subseteq  \mathcal{L}_{n}(q)$ forms a $d$-{\em dimensional} $q$-{\em algebra} if there exist
pairwise disjoint sets $X_{0}, X_{1}, \dots, X_{d} \subseteq \mathbb{F}_{q}^{n}$, all non-empty with perhaps the exception of $X_{0}$,
so that
\[\mathcal{B} = \bigg\{span\{X_{0} \cup(\cup_{i \in I}X_{i})\} \mid I \subseteq [d]\bigg\}.\]
\end{definition}

\begin{definition}\label{def5.7}
Given a positive integer d,
define $b_{q}(n, d)$ to be the maximum size of a family $\mathcal{V} \subseteq \mathcal{L}_{n}(q)$ which does
not contain a $d$-dimensional $q$-algebra.
\end{definition}

Note that forbidden $d$-dimensional $q$-algebra is a hereditary property. Theorems \ref{thm1.9} and \ref{thm5.5} together imply immediately
the next $q$-analogue of Theorem \ref{thm5.5}, where $l_{q}(\mathcal{V})$ is given in (4).

\begin{theorem}\label{thm5.8}
For $d \geq 3$ and $n \geq (2^{d} - \frac{2}{ln2})^{2}$, every family $\mathcal{V} \subseteq \mathcal{L}_{n}(q)$ containing no $d$-dimensional $q$-algebra satisfies
\[l_{q}(\mathcal{V}) \leq 2(n+1)^{1 - 2^{1-d}}.\]
\end{theorem}

Recall that $\qbinom{n}{k}$ is symmetric and has the unimodal property, we have $\qbinom{n}{k} \leq \qbinom{n}{\lceil \frac{n}{2}\rceil}$ for all $k \leq n$. Theorem \ref{thm5.8} gives the following immediate consequence which is a generalization of Theorems \ref{thm5.3} and \ref{thm5.4}.

\begin{theorem}\label{thm5.9}
For $d \geq 3$ and $n \geq (2^{d} - \frac{2}{ln2})^{2}$, every family $\mathcal{V} \subseteq \mathcal{L}_{n}(q)$ containing no $d$-dimensional $q$-algebra satisfies
\[|\mathcal{V}| \leq 2(n+1)^{1 - 2^{1-d}}\qbinom{n}{\lceil \frac{n}{2}\rceil}.\]
\end{theorem}

Gunderson et al. \cite{grs} and Johnston et al. \cite{jlm} introduced the next Ramsey-type concept: Given positive integers $n$ and $d$, define $r(d, n)$ to be the largest integer $r$ so that every $r$-coloring of $2^{[n]}$ contains a monochromatic copy of a $d$-dimensional Boolean algebra $\mathcal{B}_{d}$.
Similarly, one defines  $r_{q}(d, n)$ to be the largest integer $r$ so that every $r$-coloring of $\mathcal{L}_{n}(q)$ contains a monochromatic copy of a $d$-dimensional $q$-algebra. The following bounds for $r(d, n)$ are given in \cite{grs} and \cite{jlm}:
\[\bigg\lfloor\frac{1}{2}n^{\frac{2}{2^{d}}}\bigg\rfloor \leq r(d, n) \leq n^{\frac{d}{2^{d}-1}(1 + o(1))}.\]
By applying Theorem \ref{thm5.8}, we can derive the next lower bound for  $r_{q}(d, n)$.

\begin{theorem}\label{thm5.10}
For $d \geq 3$ and $n \geq (2^{d} - \frac{2}{ln2})^{2}$, we have
\[r_{q}(d, n) \geq \bigg\lfloor\frac{1}{2}n^{\frac{2}{2^{d}}}\bigg\rfloor.\]
\end{theorem}

\noindent {\bf Proof.}
Let $r = \bigg\lfloor\frac{1}{2}n^{\frac{2}{2^{d}}}\bigg\rfloor$.
For every $r$-coloring of $\mathcal{L}_{n}(q)$ and $1 \leq i \leq r$, let $\mathcal{V}_{i}$
be the family of subspaces in color $i$. By the definition of $l_{q}(\mathcal{V})$ (see (4)) and (24), we have
\[\sum_{1 \leq i \leq r}l_{q}(\mathcal{V}_{i}) = l_{q}(\mathcal{L}_{n}(q)) = n + 1.\]
It follows that there exists a color $1 \leq i \leq r$ such that
$l_{q}(\mathcal{V}_{i}) \geq \frac{n+1}{r} > 2(n+1)^{1-2^{1-d}}$.
By Theorem \ref{thm5.8}, $\mathcal{V}_{i}$ contains a $d$-dimensional $q$-algebra.
$\hfill \Box$

\vspace{3mm}

\section{Concluding Remarks}
We provided a general relationship for LYM inequalities between Boolean lattice and linear lattice (Theorem \ref{thm1.9}) in section 2. As applications of Theorem \ref{thm1.9}, we derived generalizations of some well-known theorems on LYM inequalities and diamond-free families of subsets of $[n]$ in section 3; the generalizations of the well-known Kleitman Theorem on matchings (Theorem \ref{thm1.13}) and its generalization (Theorem \ref{thm1.14}) in section 4; and Theorems \ref{thm5.8} and \ref{thm5.9} on forbidden $d$-dimensional $q$-algebra in section 5.

In 1975, Szemer\'edi \cite{s} proved the following famous theorem.

\begin{theorem}\label{thm6.1} (Szemer\'edi, \cite{s})
For every positive integer $k$ and every $\delta > 0$ there exists $N$ such that every
subset $A \subseteq [N]$ of size at least $\delta N$ contains an arithmetic progression of length $k$.
\end{theorem}

Szemer\'edi theorem has now three substantially different proofs and has several generalitions, including the multidimensional Szemer\'edi theorem of Furstenberg
and Katznelson \cite{fk1} (and a different proof was given by Gowers \cite{g1}).
We propose the following conjecture for a vector space variant of Szemer\'edi theorem, where we say a family $\mathcal{V}$ of subspaces of the $n$-dimensional vector space $\mathbb{F}_{q}^{n}$ contains an arithmetic progression of length $k$ if $\mathcal{V}$ contains a subfamily $\mathcal{V}'$ such that the dimensions of the
subspaces in $\mathcal{V}'$ form an arithmetic progression of length $k$.

\begin{conjecture}\label{conj6.2}
For every positive integer $k$ and every $\delta > 0$ there exists $n$ such that every
family $\mathcal{V}$ of subspaces of $\mathbb{F}_{q}^{n}$ with desity at least $\delta$ contains an arithmetic progression of length $k$.
\end{conjecture}

Recall that Sperner's theorem and its multidimensional version (Theorem \ref{thm5.3}) were applied in a new proof of Furstenberg and Katznelson's density Hales-Jewett theorem by Polymath \cite{p}.
We wonder if one could obtain a proof for Conjecture \ref{conj6.2} by following Polymath's  approach and applying Theorem \ref{thm5.9} (a generalization of Theorems \ref{thm5.3}).

\vspace{3mm}

\noindent {\bf Appendix 1}

\vspace{3mm}

\noindent {\bf Proof of Proposition 4.11.}
Recall that $\varphi_{q}(i) = \frac{|\mathcal{V}_{i}|}{\qbinom{n}{i}}$ for $0 \leq i \leq n$.
Note that since any family $\mathcal{V}$ of subspaces of $\mathbb{F}_{q}^{n}$ is a complex, $\mathbb{F}_{q}^{n} \in \mathcal{V}$ implies that
$\mathcal{V} = \mathcal{L}_{n}(q)$ and so $||\mathcal{V}||_{q} = n+1$ by (13) and (24). Hence we may assume $\varphi_{q}(n) = \frac{|\mathcal{V}_{n}|}{\qbinom{n}{n}} = 0$ for proving $q$-perfectness.
Define
\[\alpha(j) = \varphi_{q}(j) - \varphi_{q}(j+1) \mbox{ for } 0 \leq j \leq n - 1.\]
By Lemma \ref{lem4.8}, $\alpha(j) \geq 0$ for every $0 \leq j \leq n - 1$. Since $\varphi_{q}(n) = 0$, by the definition, we have
\begin{equation}
\varphi_{q}(j) = \alpha(j) + \alpha(j+1) + \cdots + \alpha(n-1).
\end{equation}
Let $\mathcal{V} \subseteq \mathcal{L}_{n}(q)$ be a complex with $||\mathcal{V}||_{q} < n+1$.
Define
\[\mathcal{V}^{+} = \{V \in \mathcal{V} \mid dim(V) \geq l\},\]
\[\mathcal{V}^{-} = \{V \in \mathcal{V} \mid dim(V) < l\}.\]
Then $\mathcal{V} = \mathcal{V}^{+} \cup \mathcal{V}^{-}$.
Multiplying both sides of (15) by $\alpha(k)$ gives
\begin{equation}
\sum_{l \leq j \leq k} \alpha(k)\qbinom{n}{j} \geq \alpha(k)(k - l + 1)\qbinom{n}{l}.
\end{equation}
By (28), we have
\begin{equation}
||\mathcal{V}^{+}||_{q} = \sum_{l \leq j \leq n-1}\varphi_{q}(j) = \sum_{l \leq j \leq n-1}\sum_{j \leq k \leq n - 1}\alpha(k)
= \sum_{l \leq k \leq n - 1}(k - l + 1)\alpha(k).
\end{equation}
It follows from (29) and (30) that
\[||\mathcal{V}^{+}||_{q}\qbinom{n}{l} \leq \sum_{l \leq k \leq n - 1}\alpha(k)\sum_{l \leq j \leq k}\qbinom{n}{j} \hspace{22mm}\]
\[= \sum_{l \leq j \leq n-1}\qbinom{n}{j}\sum_{j \leq k \leq n - 1}\alpha(k)\]
\[ = \sum_{l \leq j \leq n-1}\qbinom{n}{j}\varphi_{q}(j) \hspace{13mm}\]
\begin{equation}
= \sum_{l \leq j \leq n-1}|\mathcal{V}_{j}| = |\mathcal{V}^{+}|. \hspace{8mm}
\end{equation}
For $\mathcal{V}^{-}$,
since $\varphi_{q}(j) \leq 1$ for each $0 \leq j \leq n$, $||\mathcal{V}^{-}||_{q} = \sum_{0 \leq j < l}\varphi_{q}(j)$, and
$\qbinom{n}{j} < \qbinom{n}{l}$ for each $j < l$, we have
\[|\mathcal{V}^{-}| = \sum_{0 \leq j < l}\varphi_{q}(j)\qbinom{n}{j} \hspace{38mm}\]
\[= \sum_{0 \leq j < l}\qbinom{n}{j} - \sum_{0 \leq j < l}(1 - \varphi_{q}(j))\qbinom{n}{j}\]
\[\geq \sum_{0 \leq j < l}\qbinom{n}{j} - \sum_{0 \leq j < l}(1 - \varphi_{q}(j))\qbinom{n}{l}\]
\begin{equation}
 = \sum_{0 \leq j < l}\qbinom{n}{j} - (l - ||\mathcal{V}^{-}||_{q})\qbinom{n}{l}. \hspace{6mm}
\end{equation}
Observe that $|\mathcal{V}| = |\mathcal{V}^{+}| + |\mathcal{V}^{-}|$ and $||\mathcal{V}||_{q} = ||\mathcal{V}^{+}||_{q} + ||\mathcal{V}^{-}||_{q}$.
It follows from (31) and (32) that
\[|\mathcal{V}| \geq \sum_{0 \leq i < l}\qbinom{n}{i} + (||\mathcal{V}||_{q} - l)\qbinom{n}{l}\]
which is (25).
$\hfill \Box$

\vspace{3mm}

\section*{Declaration of Competing Interest}
The authors declare that they have no conflicts of interest to this work.

\section*{Acknowledgement}
The research is supported by the National Natural Science Foundation of China (71973103, 11861019) and
Guizhou Talent Development Project in Science and Technology (KY[2018]046), Natural Science Foundation of Guizhou ([2019]1047, [2020]1Z001,  [2021]5609).


\vspace{3mm}


\begin{thebibliography}{99}
\bibitem{abs} N. Alon, L. Babai, and H. Suzuki,
Multilinear Polynomials and Frankl-Ray-Chaudhuri-Wilson Type
Intersection Theorems, {\em J. Combin. Theory, Ser. A}, 58 (1991),
165-180.

\bibitem{b1} C. Bey,
Polynomial LYM inequalities, {\em Combinatorica}, 25 (2005), 19-38.

\bibitem{b} B. Bollob\'as,
On generalized graphs, {\em Acta Math. Acad. Sci. Hungar.}, 16 (1965),
447-452.

\bibitem{cp} A. Chowdhury and B. Patk\'os,
Shadows and intersections in vector spaces, {\em J. Combin. Theory, Ser. A}, 117 (2010),
1095-1106.

\bibitem{cdjs} E. Czabarka, A. Dutle, T. Johnston, and L. A. Sz\'ekely,
Abelian groups yield many large families for the diamond problem, {\em Eur. J. Math.}, 1 (2015),
320-328.

\bibitem{e2} K. Engel,
Sperner Theory, Volume 65 of Encyclopedia of Mathematics and its Applications,
Cambridge University Press, Cambridge, 1997

\bibitem{e} P. Erd\H{o}s,
On a lemma of Littlewood and Offord, {\em Bulletin of the American Mathematical Society}, 52 (1945), 898-902.


\bibitem{e1} P. Erd\H{o}s,
A problem on independent $r$-tuples, {\em Ann. Univ. Sci. Budapest. E\"otv\"os Sect. Math.}, 8 (1965), 83-95.

\bibitem{ekr} P. Erd\H{o}s, C. Ko, and R. Rado,
Intersection theorems for systems of finite sets, {\em Q. J. Math.
Oxford (2)}, 12 (1961), 313-320.

\bibitem{ek} P. Erd\H{o}s and D. Kleitman,
On collections of subsets containing no $4$-member Boolean algebra, {\em Proc. Amer. Math. Soc.}, 28 (1971), 87-90.

\bibitem{efkss} P. L. Erd\H{o}s, P. Frankl, D. J. Kleitman, M. E. Saks, and L. A. Sz$\acute{e}$kely,
Sharpening the LYM inequality, {\em Combinatorica}, 2 (1992), 287-293.

\bibitem{f} P. Frankl,
Extremal set theory for the binomial norm, {\em J. Combin. Theory, Ser. A}, 170 (2020), 105144.

\bibitem{fg} P. Frankl and R. Graham,
Intersection theorems for vector spaces, {\em Europ. J. Combin.}, 6 (1985), 183-187.

\bibitem{fk} P. Frankl and A. Kupavskii,
Two problems on matchings in set families - in the footsteps of Erd\H{o}s and Kleitman, {\em J. Combin. Theory, Ser. B}, 138 (2019), 286-313.

\bibitem{fw} P. Frankl and R. M. Wilson,
Intersection theorems with geometric consequences, {\em
Combinatorica}, 1 (1981), 357-368.

\bibitem{fk1} H. Furstenberg and Y. Katznelson,
An ergodic Szemer\'edi theorem for commuting transformations, {\em J. Analyse Math.}, 34 (1978), 275-291.


\bibitem{g} D. Gerbner,
The covering lemma and $q$-analogues of extremal set
theory problems, {\em Ars Mathematica Contemporanea}, 2023, doi.org/10.26493/1855-3974.2677.b7f.

\bibitem{g1} W. T. Gowers,
Hypergraph regularity and the multidimensional Szemer\'edi theorem, {\em Annals
 of Mathematics}, 166 (2007), 897-946.


\bibitem{gk} C. Greene and D. J. Kleitman,
Proof techniques in the theory of finite sets, {\em MAA Studies in Math.}, 17 (1978),
12-79.

\bibitem{gll} J. R. Griggs, W. Li, and L. Lu,
Diamond-free families, {\em J. Combin. Theory, Ser. A}, 119 (2012), 310-322.

\bibitem{gmt} D. Gr\'osz, A. Methuku, and C. Tompkins,
An upper bound on the size of diamond-free families of sets, {\em J. Combin. Theory, Ser. A}, 156 (2018), 164-194.

\bibitem{grs} D. S. Gunderson, V. R\"odl, and A. Sidorenko,
Extremal Problems for Sets Forming Boolean Algebras
and Complete Partite Hypergraphs, {\em J. Combin. Theory, Ser. A}, 88 (1999), 342-367.

\bibitem{h} W. N. Hsieh,
Intersection theorems for systems of finite vector spaces,
{\em Discrete Math.}, 12 (1975), 1-16.

\bibitem{jlm} T. Johnston, L. Lu and K. G. Milans,
Boolean algebras and Lubell functions, {\em J. Combin. Theory, Ser. A}, 136 (2015), 174-183.

\bibitem{kr} D. A. Klain and G. C. Rota,
A Continuous Analogue of Sperner's Theorem, {\em Communications on Pure and Applied Mathematics}, L (1997), 205-223.

\bibitem{k} D. J. Kleitman,
Maximal number of subsets of a finite set no $k$ of which are pairwise disjoint, {\em J. Combin. Theory}, 5 (1968), 157-163.


\bibitem{kmy} L. Kramer, R. R. Martin, and M. Young,
On diamond-free subposets of the Boolean lattice, {\em J. Combin. Theory, Ser. A}, 120 (2013), 545-560.

\bibitem{lm} L. Lu and K. G. Milans,
Set families with forbidden subposets, {\em J. Combin. Theory, Ser. A}, 136 (2015), 126-142.

\bibitem{l2} D. Lubell,
A short proof of Sperner's lemma, {\em J. Combin. Theory}, 1 (1966), 299.

\bibitem{mt} D. Malec and C. Tompkins,
Localized versions of extremal problems, {\em European Journal of Combinatorics}, 112 (2023), 103715.

\bibitem{m} A. M$\grave{e}$roueh,
A LYM inequality for induced posets, {\em J. Combin. Theory, Ser. A}, 155 (2018), 398-417.

\bibitem{m1} L. D. Meshalkin,
Generalization of Sperner's theorem on the number of subsets of a finite set, {\em Theory Probab. Appl.}, 8 (1963), 203-204.

\bibitem{p} D. H. J. Polymath,
A new proof of the density Hales-Jewett theorem, {\em Annals of Mathematics}, 175 (2012), 1283-1327.

\bibitem{qr} J. Qian and D. K. Ray-Chaudhuri,
Extremal case of Frankl-Ray-Chaudhuri-Wilson inequality, {\em J. Statist. Plann. Inference}, 95 (2001), 293-306.

\bibitem{rw} D. K. Ray-Chaudhuri and R. M. Wilson,
 On $t$-designs, {\em Osaca J. Math.}, 12 (1975), 737-744.

\bibitem{s2} W. Samotij,
Subsets of Posets Minimising the Number of Chains, {\em Transactions of the American Mathematical Society}, 371 (2019), 7259-7274.

\bibitem{s} E. Szemer\'edi,
On sets of integers containing no k elements in arithmetic progression, {\em Acta Arithmetica}, 27 (1975), 199-245.

\bibitem{y} K. Yamamoto,
Logarithmic order of free distributive lattice, {\em J. Math. Soc. Japan}, 6 (1954), 343-353.




\end{thebibliography}
\end{document}